\documentclass[11pt]{article}
\usepackage{mathrsfs}
\usepackage{amssymb}
\usepackage{amsmath}
\usepackage[all]{xy}
\def\Im{\mathop{\rm Im}\nolimits}
\def\Ker{\mathop{\rm Ker}\nolimits}
\def\Coker{\mathop{\rm Coker}\nolimits}
\def\Tr{\mathop{\rm Tr}\nolimits}
\def\mod{\mathop{\rm mod}\nolimits}
\def\Mod{\mathop{\rm Mod}\nolimits}
\def\fd{\mathop{\rm fd}\nolimits}
\def\id{\mathop{\rm id}\nolimits}
\def\pd{\mathop{\rm pd}\nolimits}
\def\max{\mathop{\rm max}\nolimits}

\def\sup{\mathop{\rm sup}\nolimits}

\def\length{\mathop{\rm length}\nolimits}
\def\gldim{\mathop{\rm gl.dim}\nolimits}
\def\findim{\mathop{\rm fin.dim}\nolimits}
\def\Rapp{\mathop{\rm Rapp}\nolimits}
\def\Lapp{\mathop{\rm Lapp}\nolimits}
\def\G{\mathop{\mathscr{G}}\nolimits}
\def\CoG{\mathop{{\rm Co}\mathscr{G}}\nolimits}
\def\Findim{\mathop{\rm Fin.dim}\nolimits}
\def\Hom{\mathop{\rm Hom}\nolimits}
\def\Ext{\mathop{\rm Ext}\nolimits}
\def\sup{\mathop{\rm sup}\nolimits}
\def\lim{\mathop{\underrightarrow{\rm lim}}\nolimits}
\def\colim{\mathop{\underleftarrow{\rm lim}}\nolimits}

\title{\Large \bf On Auslander-Type Conditions of Modules
}
\author{Zhaoyong Huang\thanks{E-mail address: huangzy@nju.edu.cn}\\
{\it \footnotesize Department of Mathematics, Nanjing University,
Nanjing 210093, Jiangsu Province, P. R. China}}
\date{ }
\begin{document}
\baselineskip=18pt \maketitle

\begin{abstract}
We prove that for a left and right Noetherian ring $R$, $_RR$
satisfies the Auslander condition if and only if so does every flat
left $R$-module, if and only if the injective dimension of the $i$th
term in a minimal flat resolution of any injective left $R$-module
is at most $i-1$ for any $i \geq 1$, if and only if the flat (resp.
injective) dimension of the $i$th term in a minimal injective coresolution (resp.
flat resolution) of any left $R$-module $M$ is at most the flat
(resp. injective) dimension of $M$ plus $i-1$ for any $i \geq 1$, if
and only if the flat (resp. injective) dimension of the injective
envelope (resp. flat cover) of any left $R$-module $M$ is at most
the flat (resp. injective) dimension of $M$, and if and only if any
of the opposite versions of the above conditions hold true.
Furthermore, we prove that for an Artinian algebra $R$ satisfying
the Auslander condition, $R$ is Gorenstein if and only if the
subcategory consisting of finitely generated modules satisfying the
Auslander condition is contravariantly finite.
As applications, we get some equivalent
characterizations of Auslander-Gorenstein rings and
Auslander-regular rings.
\end{abstract}

\vspace{0.5cm}

\centerline{\large \bf 1. Introduction}

\vspace{0.2cm}

It is well known that commutative Gorenstein rings are fundamental
and important research objects in commutative algebra and algebraic
geometry. Bass proved in [B2] that a commutative Noetherian ring $R$
is a Gorenstein ring (that is, the self-injective dimension of $R$
is finite) if and only if the flat dimension of the $i$th term in a
minimal injective coresolution of $R$ as an $R$-module is at most
$i-1$ for any $i \geq 1$. In non-commutative case, Auslander proved
that this condition is left-right symmetric ([FGR, Theorem 3.7]). In
this case, $R$ is said to satisfy the {\it Auslander condition}.
Motivated by this philosophy, Huang and Iyama introduced the notion
of Auslander-type conditions of rings as follows. For any $m,n\geq
0$, a left and right Noetherian ring is said to be $G_n(m)$ if the
flat dimension of the $i$th term in a minimal injective coresolution
of $R_R$ is at most $m+i-1$ for any $1\leq i \leq n$. Auslander-type
conditions are non-commutative analogs of commutative Gorenstein
rings. Such conditions play a crucial role in homological algebra,
representation theory of algebras and non-commutative algebraic
geometry ([AR3], [AR4], [Bj], [EHIS], [FGR], [H1], [HI], [IS], [I1],
[I2], [I3], [I4], [M], [Ro], [S], [W], and so on). In particular, by
constructing an injective coresolution of the last term in an exact
sequence of finite length from that of the other terms, Miyachi
obtained in [M] an equivalent characterization of the Auslander
condition in terms of the relation between the flat dimensions of
any module and its injective envelope. Then he got some properties
of Auslander-Gorenstein rings and Auslander-regular rings.

Note that a commutative Noetherian ring satisfies the Auslander
condition if and only if it is Gorenstein ([B2]). Auslander and
Reiten conjectured in [AR3] that an Artinian algebra satisfying the
Auslander condition is Gorenstein. This conjecture is situated
between the well known Nakayama conjecture and the finitistic
dimension conjecture. For an Artinian algebra $R$, the Nakayama
conjecture states that $R$ is selfinjective if all terms in a
minimal injective coresolution of $_RR$ are projective; and the
finitistic dimension conjecture states that the supremum of the
projective dimensions of all finitely generated left $R$-modules
with finite projective dimension is finite. All of these conjectures
remains still open.

Based on these mentioned above, in this paper we will introduce
modules satisfying Auslander-type conditions and study the
homological properties of such modules. By using the obtained
properties, we get some equivalent characterizations of rings
satisfying the Auslander condition, Auslander-Gorenstein rings and
Auslander-regular rings respectively. Then we study when an Artinian
algebra satisfying the Auslander condition is Gorenstein. This paper
is organized as follows.

In Section 2, we give some terminology and some preliminary results.

In Section 3, by using some techniques of direct limits and
transfinite induction, we prove the following

\vspace{0.2cm}

{\bf Theorem 1.1.} {\it Let $R$ be a left Noetherian ring and
$n,k\geq 0$, and let $\{M_i\}_{i\in I}$ be a family of left
$R$-modules and $M=\underset{i\in I}\lim M_i$, where $I$ is a
directed index set. If the flat dimension of the $(n+1)$st term in a
minimal injective coresolution of $M_i$ is at most $k$ for any $i\in
I$, then the flat dimension of the $(n+1)$st term in a minimal
injective coresolution of $M$ is also at most $k$.}

\vspace{0.2cm}

For any $m,n\geq 0$, we introduce in Section 4 the notion of modules
satisfying the Auslander-type conditions $G_n(m)$; in particular, a
left $R$-module $M$ for any ring $R$ is said to satisfy the {\it
Auslander condition} if the flat dimension of the $i$th term in a
minimal injective coresolution of $_RM$ is at most $i-1$ for any
$i\geq 1$. By using Theorem 1.1 and the constructions of (co)proper
(co)resolutions of modules in [H2], we will investigate the
homological behavior of modules satisfying Auslander-type conditions
in terms of the relation between the flat (resp. injective)
dimensions of modules and their injective envelopes (resp. flat
covers). We prove the following

\vspace{0.2cm}

{\bf Theorem 1.2.} {\it Let $R$ be a left and right Noetherian ring.
Then the following statements are equivalent.

(1) $_RR$ satisfies the Auslander condition.

(2) Every flat left $R$-module satisfies the Auslander condition.

(3) The flat dimension of the $i$th term in a minimal injective
coresolution of any left $R$-module $M$ is at most the flat dimension
of $M$ plus $i-1$ for any $i \geq 1$.

(4) The flat dimension of the injective envelope of any left
$R$-module $M$ is at most the flat dimension of $M$.

(5) The injective dimension of the $i$th term in a minimal flat
resolution of any injective left $R$-module is at most $i-1$ for
any $i \geq 1$.

(6) The injective dimension of the $i$th term in a minimal flat
resolution of any left $R$-module $M$ is at most the injective
dimension of $M$ plus $i-1$ for any $i \geq 1$.

(7) The injective dimension of the flat cover of any left $R$-module
$M$ is at most the injective dimension of $M$.

$(i)^{op}$ The opposite version of $(i)$ ($1\leq i \leq 7$).}

\vspace{0.2cm}

As applications of this theorem, we obtain some equivalent
characterizations of Auslander-Gorenstein rings and
Auslander-regular rings, respectively.

In Section 5, we first obtain the approximation presentations of a
given module relative to the subcategory of modules satisfying the
Auslander condition and that of modules with finite injective
dimension respectively. Then we establish the connection between the
Auslander and Reiten conjecture mentioned above with the
contravariant finiteness of some certain subcategories as follows.

\vspace{0.2cm}

{\bf Theorem 1.3.} {\it Let $R$ be an Artinian algebra satisfying
the Auslander condition. Then the following statements are
equivalent.

(1) $R$ is Gorenstein.

(2) The subcategory consisting of finitely generated modules satisfying the
Auslander condition is contravariantly finite.

(3) The subcategory consisting of finitely generated modules which
are $n$-syzygy for any $n\geq 1$ is contravariantly finite.

(4) The subcategory consisting of finitely generated modules which
are $n$-torsionfree for any $n\geq 1$ is contravariantly finite.}


\vspace{0.2cm}

As a consequence, we get that an Artinian algebra is
Auslander-regular if and only if the subcategory consisting of
projective modules and that consisting of modules satisfying the
Auslander condition coincide.

\vspace{0.3cm}

\centerline{\large \bf 2. Preliminaries}

\vspace{0.2cm}

Throughout this paper, $R$ is an associative ring with identity,
$\Mod R$ is the category of left $R$-modules and $\mod R$ is the
category of finitely generated left $R$-modules. We use $\gldim R$
to denote the global dimension of $R$. In this section, we give some
terminology and some preliminary results.

\vspace{0.2cm}

{\bf Definition 2.1.} ([E]) Let $\mathscr{C}\subseteq\mathscr{D}$ be
full subcategories of $\Mod R$. The homomorphism $f: C\to D$ in $\Mod R$
with $C\in\mathscr{C}$ and $D\in\mathscr{D}$ is said to be a
{\it $\mathscr{C}$-precover} of $D$ if for any homomorphism $g: C^{'} \to D$
in $\Mod R$ with $C^{'}\in\mathscr{C}$, there exists a homomorphism
$h: C^{'}\to C$ such that the following diagram commutes:
$$\xymatrix{ & C^{'} \ar[d]^{g} \ar@{-->}[ld]_{h}\\
C \ar[r]^{f} & D}$$
The homomorphism $f: C\to D$ is said to be {\it
right minimal} if an endomorphism $h: C\to C$ is an automorphism
whenever $f=fh$. A $\mathscr{C}$-precover $f: C\to D$ is called a
{\it $\mathscr{C}$-cover} if $f$ is right minimal. Dually, the
notions of a {\it $\mathscr{C}$-preenvelope}, a {\it left minimal
homomorphism} and {\it a $\mathscr{C}$-envelope} are defined. Following
Auslander and Reiten's terminology in [AR1], for a module over an Artinian
algebra, a $\mathscr{C}$-(pre)cover and a $\mathscr{C}$-(pre)envelope are
called a {\it (minimal) right $\mathscr{C}$-approximation} and a {\it (minimal)
left $\mathscr{C}$-approximation}, respectively. If each module in $\mathscr{D}$
has a right (resp. left) $\mathscr{C}$-approximation, then $\mathscr{C}$
is called {\it contravariantly finite} (resp. {\it covariantly finite})
in $\mathscr{D}$.

\vspace{0.2cm}







{\bf Lemma 2.2.} {\it ([X, Theorem 1.2.9]) Let $\mathscr{C}$ be a
full subcategory of $\Mod R$ closed under direct products. If $f_i:
C_i \to M_i$ is a $\mathscr{C}$-precover of $M_i$ in $\Mod R$ for
any $i \in I$, where $I$ is an index set, then $\prod _{i \in I}f_i:
\prod _{i \in I}C_i\to \prod _{i \in I}M_i$ is a
$\mathscr{C}$-precover of $\prod _{i \in I}M_i$.}

\vspace{0.2cm}

We use $\mathscr{F}^0(\Mod R)$ and $\mathscr{I}^0(\Mod R)$ to denote the
subcategories of $\Mod R$ consisting of flat modules and injective
modules, respectively. Recall that an $\mathscr{F}^0(\Mod R)$-(pre)cover
and an $\mathscr{I}^0(\Mod R)$-(pre)envelope are called a {\it flat
(pre)cover} and an {\it injective (pre)envelope}, respectively.

Bican, El Bashir and Enochs proved in [BEE, Theorem 3] that every
$R$-module has a flat cover. For an $R$-module $M$, we call an exact
sequence $\cdots \to F_i\buildrel {{\pi}_i} \over \longrightarrow
\cdots \buildrel {{\pi}_2} \over \longrightarrow F_1 \buildrel
{{\pi}_1} \over\longrightarrow F_0 \buildrel {{\pi}_0}
\over\longrightarrow M \to 0$ a {\it proper flat resolution} of $M$
if ${{\pi}_i}: F_i \to \Im {{\pi}_{i}}$ is a flat precover of $\Im
{{\pi}_{i}}$ for any $i \geq 0$. Furthermore, we call the following
exact sequence:
$$\cdots \to F_i(M)\buildrel {{\pi}_i(M)} \over \longrightarrow \cdots
\buildrel {{\pi}_2(M)} \over \longrightarrow F_1(M) \buildrel
{{\pi}_1(M)} \over\longrightarrow F_0(M) \buildrel {{\pi}_0(M)}
\over\longrightarrow M \to 0$$ a {\it minimal flat resolution} of
$M$, where ${{\pi}_i}(M): F_i(M) \to \Im {{\pi}_{i}}(M)$ is a flat
cover of $\Im {{\pi}_{i}}(M)$ for any $i \geq 0$. It is easy to
verify that the flat dimension of $M$ is at most $n$ if and only if
$F_{n+1}(M)=0$. In addition, we use $$0 \to M \to E^0(M) \to E^1(M)
\to \cdots \to E^i(M) \to \cdots$$ to denote a minimal injective
coresolution of $M$.

We denote by $(-)^+=\Hom _{\mathbb{Z}}(-, \mathbb{Q}/\mathbb{Z})$,
where $\mathbb{Z}$ is the additive group of integers and
$\mathbb{Q}$ is the additive group of rational numbers.

\vspace{0.2cm}

{\bf Lemma 2.3.} {\it ([EH, Theorem 3.7]) The following statements
are equivalent.

(1) $R$ is a left Noetherian ring.

(2) A monomorphism $f: A \rightarrowtail E$ in $\Mod R$ is an
injective preenvelope of $A$ if and only if $f^+: E^+
\twoheadrightarrow A^+$ is a flat precover of $A^+$ in $\Mod
R^{op}$.}

\vspace{0.2cm}

Let $M\in \Mod R$. We use $\fd _RM$, $\pd _RM$ and $\id _RM$ to
denote the flat, projective and injective dimensions of $M$,
respectively.

\vspace{0.2cm}

{\bf Lemma 2.4.} {\it (1) ([F, Theorem 2.1]) For any $M\in \Mod R$,
$\fd _RM=\id _{R^{op}}M^+$.

(2) ([F, Theorem 2.2]) If $R$ is a right Noetherian ring, then $\fd
_RN^+=\id _{R^{op}}N$ for any $N\in \Mod R^{op}$.}

\vspace{0.2cm}

Recall that $\Findim R=\sup\{\pd _RM\ |\ M\in \Mod R$ with $\pd
_RM<\infty\}$. Observe that the first assertion in the following
result was proved by Bass in [B1, Corollary 5.5] when $R$ is a
commutative Noetherian ring.

\vspace{0.2cm}

{\bf Lemma 2.5.} {\it (1) For a left Noetherian ring $R$, $\id
_RR\geq \sup\{\fd _RM\ |\ M\in \Mod R$ with $\fd _RM<\infty\}$.

(2) For a left and right Noetherian ring $R$, $\id _RR\geq \sup\{\id
_{R^{op}}N\ |\ N\in \Mod R^{op}$ with $\id _{R^{op}}N<\infty\}$.}

\vspace{0.2cm}

{\it Proof.} (1) Without loss of generality, assume that $\id
_RR=n<\infty$. Then $\Findim R\leq n$ by [B1, Proposition 4.3]. It
follows from [J1, Proposition 6] that the projective dimension of any
flat left $R$-module is finite. So, if $M\in \Mod R$ with $\fd
_RM<\infty$, then $\pd _RM<\infty$ and $\pd _RM\leq n$. Thus we have
$\fd _RM(\leq \pd _RM)\leq n$.

(2) By [B1, Proposition 4.1], we have $\sup\{\fd _RM\ |\ M\in \Mod
R$ with $\fd _RM<\infty\}=\sup\{\id _{R^{op}}N\ |\ N\in \Mod R^{op}$
with $\id _{R^{op}}N<\infty\}$. So the assertion follows from (1).
\hfill$\square$

\vspace{0.5cm}

\centerline{\large\bf 3. Flat dimension of E$^n$ of direct limits}

\vspace{0.2cm}

In this section, $R$ is a left Noetherian ring. The aim of
this section is to prove the following

\vspace{0.2cm}

{\bf Theorem 3.1.} {\it Let $n,k\geq 0$ and let $\{M_i\}_{i\in I}$
be a family of left $R$-modules, where $I$ is a directed index set.
If $M=\underset{i\in I}\lim M_i$ and $\fd _RE^n(M_i)\leq k$ for any
$i\in I$, then $\fd _RE^n(M)\leq k$.}

\vspace{0.2cm}

By [R, Theorem 5.40], every flat left $R$-module is a direct limit
(over a directed index set) of finitely generated free left
$R$-modules. So by Theorem 3.1, we have the following

\vspace{0.2cm}

{\bf Corollary 3.2.} {\it $\fd _RE^n(_RR)=\sup\{\fd _RE^n(F)\ |\
F\in\Mod R$ is flat$\}$ for any $n\geq 0$.}

\vspace{0.2cm}

Before giving the proof of Theorem 3.1, we need some preliminaries.

\vspace{0.2cm}

{\bf Definition 3.3.} ([J2]) Let $\beta$ be an ordinal number. A set
$S$ is called a {\it continuous union} of a family of subsets
indexed by ordinals $\alpha$ with $\alpha <\beta$ if for each such
$\alpha$ we have a subset $S_{\alpha}\subset S$ such that if $\alpha
\leq \alpha^{'}$ then $S_{\alpha}\subset S_{\alpha^{'}}$, and such
that if $\gamma < \beta$ is a limit ordinal then
$S_{\gamma}=\bigcup_{\alpha<\gamma} S_{\alpha}$.

\vspace{0.2cm}

A main tool in our proof will be the next result.

\vspace{0.2cm}

{\bf Lemma 3.4.} ([J2, Lemma 1.4]) {\it If $I$ is an infinite
directed index set, then for some ordinal $\beta$, $I$ can be
written as a continuous union $I=\bigcup_{\alpha <\beta}I_{\alpha}$,
where each $I_{\alpha}$ is a directed index set with the order
induced by that of $I$ and where $|I_{\alpha} | < |I|$ for each
$\alpha < \beta$.}

\vspace{0.2 cm}

This result will be useful since it will allow us to rewrite a
direct limit as a well-ordered direct limit. So if $M=\underset{i\in
I}\lim M_i$ with $I$ infinite, then write $I=\bigcup_{\alpha <\beta}I_{\alpha}$
as above, and put $M_{\alpha}=\underset{i\in I_{\alpha}}\lim M_i$.
Hence if $\alpha \leq \alpha^{'} <\beta$, since
$I_{\alpha} \subset I_{{\alpha}^{'}}$ we have an obvious map
$M_{\alpha}\rightarrow M_{{\alpha}^{'}}$. These maps then give us a
directed system $\{M_{\alpha}\}_{\alpha <\beta}$. Clearly
then $\underset{\alpha <\beta}\lim M_{\alpha}=\underset{i\in I}\lim
M_{i}$.

\vspace{0.2cm}

{\bf Proposition 3.5.} {\it Let $\beta $ be an ordinal number and
$\{M_{\alpha}\}$ a directed system of modules (indexed by $\alpha
<\beta$). If
\begin{center}
{$ \zeta_{\alpha}=:0\rightarrow M_{\alpha}\rightarrow
E^0(M_{\alpha})\rightarrow E^1(M_{\alpha}) \rightarrow \cdots $}
\end{center}
is a  minimal injective coresolution of $M_{\alpha}$ for each
$\alpha$, then these exact sequences $\zeta_{\alpha}$ are the
members of a directed system indexed by $\alpha <\beta$ in such a
way that if $\alpha \leq \alpha^{'}< \beta $ the map from the sequence
indexed by $\alpha$ into that indexed by $\alpha^{'}$ agrees with the
original map $M_{\alpha}\rightarrow M_{{\alpha}^{'}}$.}

\vspace{0.22cm}

{\it Proof.} Given an $\alpha+1 <\beta$ we can form a commutative
diagram: $$\xymatrix{0 \ar[r] & M_{\alpha} \ar[r] \ar[d]
& E^0(M_{\alpha}) \ar[r] \ar[d] & E^1(M_{\alpha}) \ar[r] \ar[d] & \cdots\\
0 \ar[r] & M_{\alpha+1} \ar[r] & E^0(M_{\alpha+1}) \ar[r] &
E^1(M_{\alpha+1}) \ar[r] & \cdots}$$ Using this observation we can
successively get maps $\zeta_0\rightarrow \zeta_1$,
$\zeta_1\rightarrow \zeta_2$, $\cdots$. So composing we get maps
$\zeta_m\rightarrow \zeta_n$ whenever $m\leq n$. Since $R$ is
left Noetherian, any direct limit of injective left $R$-modules
is injective by [B1, Theorem 1.1]. So $\lim \zeta_n$ is in fact an
injective coresolution of $\lim M_n$. We have a map $\lim M_n\rightarrow
M_{\omega}$ given by the maps $M_n\rightarrow M_\omega$ (where
$\omega$ is the least infinite ordinal). Then the above shows that
this in turn gives a map $\lim \zeta_n \rightarrow \zeta_{\omega}$.
So these maps give maps $\zeta_n \rightarrow \zeta_{\omega}$ for any
$n\geq 0$. Continuing this procedure we get the desired system. \hfill $\square$

\vspace{0.2cm}

Note that this result gives that if $\zeta$ is an injective
coresolution of $M$, then $\zeta\cong \underset{\alpha <\beta}\lim
\zeta_{\alpha}$. In particular, this gives that $E^n(M)\cong
\underset{\alpha <\beta}\lim E^n (M_{\alpha})$. This then gives that
if $\fd _RE^n(M_{\alpha})\leq k$ for each $\alpha$ then $\fd
_RE^n(M)\leq k$. In other words, Theorem 3.1 holds true when our
direct system is over the well-ordered index set of $\alpha <\beta$
for some ordinal $\beta$.

\vspace{0.2cm}

{\it Proof of Theorem 3.1.} We proceed by transfinite induction on
$|I|$. So to begin the induction we suppose that $|I|=\aleph_0$ (the
first infinite cardinal number). Then $I$ is countable, so we
suppose $I=\{i_n |n\in \mathbb{N}\}$ with $\mathbb{N}$ the set of
non-negative integers. We construct a sequence $j_0,j_1, j_2, \cdots
$ of elements in $I$ by letting $j_0=i_0$. Then we choose $j_1 $ so
that $j_1\geq j_0, i_1$. So in general we choose $j_n$ so that
$j_n\geq j_{n-1},i_n$. Then let $J=\{j_n |n \in \mathbb{N}\}$. We
have that $J$ is well-ordered and is clearly a confinal subset of
$I$. Hence $M=\underset{i\in I}\lim M_i =\underset{j\in J}\lim M_j$.
Since $J$ is well-ordered, $E^n(M)=\underset{j\in J}\lim E^n (M_j)$.
So the assumption that $\fd _RE^n (M_j)\leq k$ for each $j$ gives
that $\fd _RE^n(M)\leq k$.

Now we make the induction hypothesis and assume $|I|> \aleph_0$. We
appeal to Lemma 3.4 and write $I=\bigcup_{\alpha<\beta}I_{\alpha}$
as in that lemma. Then $M=\underset{\alpha<\beta}\lim M_{\alpha}$.
We have $M_{\alpha}$ is the limit over $I_{\alpha}$. But
$|I_{\alpha}|<|I|$, so the assertion holds true for direct limits
over $I_{\alpha}$ by the induction hypothesis. This means that we
have $\fd _RM_{\alpha}\leq k$ for each $\alpha$. Because the system
$\{M_{\alpha}\}_{\alpha<\beta}$ is over a well-ordered index set of
indices, we get that $\fd _RE^n (M_{\alpha})\leq k$ for each
$\alpha$ gives the assertion that $\fd _RE^n(M)\leq k$.
\hfill$\square$

\vspace{0.2cm}

{\it Remark 3.6.} The same techniques show that if for a given
$n\geq 0$ we let
$$0\rightarrow M_{\alpha}\rightarrow E^0(M_{\alpha})\rightarrow
E^1(M_{\alpha}) \rightarrow \cdots \rightarrow E^{n-1}(M_{\alpha}
)\rightarrow C^n (M_{\alpha}) \rightarrow 0$$
be a partial minimal injective coresolution of $M_{\alpha}$ for each
$\alpha$. If $\fd _RC^n(M_{\alpha})\leq k$ for each $\alpha$, then
we get that $\fd _RC^n(M)\leq k$, where
$$0\rightarrow M\rightarrow E^0(M)\rightarrow E^1(M)\rightarrow
\cdots \rightarrow E^{n-1}(M)\rightarrow C^n(M)\rightarrow 0$$
is a partial minimal injective coresolution of $M$.

\vspace{0.5cm}

\centerline{\large\bf 4. Modules satisfying the Auslander-type
conditions}

\vspace{0.2cm}

As a generalization of rings satisfying the Auslander condition,
Huang and Iyama introduced in [HI] the notion of rings satisfying
Auslander-type conditions. Now we introduce the notion of modules
satisfying the Auslander-type conditions as follows.

\vspace{0.2cm}

{\bf Definition 4.1.} Let $M\in \Mod R$ and let $m$ and $n$ be
non-negative integers. $M$ is said to be {\it $G_n(m)$} if $\fd
_RE^i(M)\leq m+i$ for any $0 \leq i \leq n-1$, and $M$ is said to be
{\it $G_{\infty}(m)$} if it is {\it $G_n(m)$} for all $n$.

\vspace{0.2cm}

{\it Remark 4.2.} Let $R$ be a left and right Noetherian ring. Then
we have

(1) $_RR$ is $G_n(m)$ if and only if $R$ is $G_n(m)^{op}$ in the
sense of Huang and Iyama in [HI].

(2) Recall from [FGR] that $R$ is called {\it Auslander's
$n$-Gorenstein} if $\fd _RE^i(_RR)\leq i$ for any $0 \leq i \leq
n-1$, and $R$ is said to satisfy the {\it Auslander condition} if it
is Auslander's $n$-Gorenstein for all $n$. So $R$ is Auslander's
$n$-Gorenstein if and only if $_RR$ is $G_n(0)$. Note that the
notion of Auslander's $n$-Gorenstein rings (and hence that of the
Auslander condition) is left-right symmetric ([FGR, Theorem 3.7]).
So $R$ satisfies the Auslander condition if and only if both $_RR$
and $R_R$ are $G_{\infty}(0)$. However, in general, the notion of
$R$ being $G_n(m)$ is not left-right symmetric when $m\geq 1$ ([AR4]
or [HI]).

\vspace{0.2cm}

The aim of this section is to study the homological behavior of
modules (especially, $_RR$) satisfying certain Auslander-type
conditions. We begin with the following

\vspace{0.2cm}

{\bf Lemma 4.3.} {\it (1) $\fd _RE^0(M)\leq \fd _RM$ for any $M\in
\Mod R$ if and only if $\fd _RE^i(M)\leq \fd _RM+i$ for any $M\in
\Mod R$ and $i\geq 0$.

(2) $\id _{R^{op}}F_0(N)\leq \id _{R^{op}}N$ for any $N\in \Mod
R^{op}$ if and only if $\id _{R^{op}}F_i(N)\leq \id _{R^{op}}N+i$
for any $N\in \Mod R^{op}$ and $i\geq 0$.}

\vspace{0.2cm}

{\it Proof.} (1) The necessity is trivial. We next prove the
sufficiency. Without loss of generality, assume that $M\in \Mod R$
with $\fd _RM=s<\infty$. In a minimal injective coresolution
$$0\to M \to E^0(M)\to E^1(M) \to \cdots \to E^i(M)\to \cdots$$
of $M$ in $\Mod R$, putting $K_{i+1}=\Im(E^i(M)\to E^{i+1}(M))$ for
any $i \geq 0$. By assumption, $\fd _RE^0(M)\leq \fd _RM=s$. So $\fd
_RK_1\leq s+1$ and hence $\fd _RE^1(M)=\fd _RE^0(K_1)\leq \fd
_RK_1\leq s+1$ again by assumption. Then $\fd _RK_2\leq s+2$.
Continuing this process, we get that $\fd _RE^i(M)\leq s+i$ for any
$i \geq 0$.

(2) It is dual to (1). \hfill $\square$

\vspace{0.2cm}

The following lemma plays an important role in the proof of the main
result of this section.

\vspace{0.2cm}

{\bf Lemma 4.4.} {\it For a left Noetherian ring $R$, $\id
_{R^{op}}F_i(E)\leq \fd _RE^i(_RR)$ for any injective right
$R$-module $E$ and $i\geq 0$.}

\vspace{0.2cm}

{\it Proof.} By Lemma 2.3, we have that
$$\cdots \to [E^i(_RR)]^+ \buildrel {{\pi}_i} \over \longrightarrow \cdots
\buildrel {{\pi}_2} \over \longrightarrow [E^1(_RR)]^+ \buildrel
{{\pi}_1} \over \longrightarrow [E^0(_RR)]^+ \buildrel {{\pi}_0}
\over \longrightarrow (_RR)^+ \to 0$$ is a proper flat resolution of
$(_RR)^+$ in $\Mod R^{op}$.

Let $E$ be an injective right $R$-module. Because $(_RR)^+$ is an
injective cogenerator for $\Mod R^{op}$, $E$ is isomorphic to a
direct summand of $[(_RR)^+]^I$ for some index set $I$. Because the
subcategory of $\Mod R^{op}$ consisting of flat modules is closed
under direct products by [C, Theorem 2.1], ${\pi _i}^I:
([E^i(_RR)]^+)^I\to (\Im {\pi}_i)^I$ is a flat precover of $(\Im
{\pi}_i)^I$ for any $i \geq 0$ by Lemma 2.2. Note that $F_i(E)$ is
isomorphic to a direct summand of $([E^i(_RR)]^+)^I$ for any $i \geq
0$. So by Lemma 2.4(1), we have that $\id _{R^{op}}F_i(E)\leq \id
_{R^{op}}([E^i(_RR)]^+)^I=\id _{R^{op}}[E^i(_RR)]^+=\fd _RE^i(_RR)$
for any $i \geq 0$. \hfill$\square$

\vspace{0.2cm}

As a consequence of Lemma 4.4 and [H2, Corollary 3.3], we get the
following

\vspace{0.2cm}

{\bf Proposition 4.5.} {\it Let $R$ be a left Noetherian ring. If
$_RR$ is $G_{\infty}(m)$ for a non-negative integer $m$, then $\id
_{R^{op}}F_i(N)\leq \id _{R^{op}}N+m+i$ for any $N\in\Mod R^{op}$
and $i \geq 0$.}

\vspace{0.2cm}

{\it Proof.} Without loss of generality, assume that $\id
_{R^{op}}N=s<\infty$. We will proceed by induction on $s$. Assume
that $_RR$ is $G_{\infty}(m)$, that is, $\fd _RE^i(_RR)\leq m+i$ for
any $i \geq 0$. If $s=0$, then the assertion follows from Lemma 4.4.

Now suppose $s\geq 1$. Then we have an exact sequence:
$$0\to N \to E^0(N) \to N_1 \to 0$$ in $\Mod R^{op}$ with $\id
_{R^{op}}N_1=s-1$. By the induction hypothesis, we have that $\id
_{R^{op}}F_i(N_1)\leq (s-1)+m+i$ and $\id _{R^{op}}F_i(E^0(N))\leq
m+i$ for any $i \geq 0$. By [H2, Remark 2.3(3) and Corollary 3.3],
we have that
$$\cdots \to F_{i+1}(N_1)\bigoplus F_i(E^0(N)) \to \cdots \to
F_2(N_1)\bigoplus F_1(E^0(N))\to F_0 \to N\to 0$$ is a strongly proper flat
resolution of $N$, and
$$0 \to F_0 \to F_1(N_1)\bigoplus F_0(E^0(N)) \to F_0(N_1)\to 0$$ is
exact. So $\id _{R^{op}}F_0\leq s+m$, and $\id
_{R^{op}}F_{i+1}(N_1)\bigoplus F_i(E^0(N))\leq s+m+i$ for any $i
\geq 1$. Notice that $F_0(N)$ is isomorphic to a direct summand of
$F_0$ and $F_i(N)$ is isomorphic to a direct summand of
$F_{i+1}(N_1)\bigoplus F_i(E^0(N))$ for any $i \geq 1$, thus we have
$\id _{R^{op}}F_i(N)\leq s+m+i$ for any $i \geq 0$. \hfill$\square$

\vspace{0.2cm}

Similarly, we get the following

\vspace{0.2cm}

{\bf Proposition 4.6.} {\it For a non-negative integer $m$, $\id
_{R^{op}}F_i(E)\leq m+i$ for any injective right $R$-module $E$ and
$i \geq 0$ if and only if $\id _{R^{op}}F_i(N)\leq \id
_{R^{op}}N+m+i$ for any $N\in\Mod R^{op}$ and $i \geq 0$.}

\vspace{0.2cm}

The following result can be regarded as a dual version of
Proposition 4.6.

\vspace{0.2cm}

{\bf Proposition 4.7.} {\it For a non-negative integer $m$, any flat
left $R$-module is $G_{\infty}(m)$ if and only if $\fd _RE^i(M)\leq
\fd _RM+m+i$ for any left $R$-module $M$ and $i \geq 0$.}

\vspace{0.2cm}

{\it Proof.} The sufficiency is trivial. We next prove the
necessity. Without loss of generality, assume that $\fd
_RM=s<\infty$. We will proceed by induction on $s$.

If $s=0$, then the assertion follows from the assumption. Now
suppose $s \geq 1$. Then we have an exact sequence:
$$0\to M_1 \to F_0(M) \to M \to 0$$ in $\Mod R$ with $\fd_RM_1=
s-1$. So by the induction hypothesis, we have that $\fd
_RE^i(M_1)\leq (s-1)+m+i$ and $\fd _RE^i(F_0(M)) \leq m+i$ for any
$i \geq 0$.

By [H2, Corollary 3.5], we have that
$$0\to M \to I^0 \to E^1(F_0(M))\bigoplus E^2(M_1)\to \cdots
\to E^i(F_0(M))\bigoplus E^{i+1}(M_1)\to \cdots$$ is an injective
coresolution of $M$, and
$$0\to E^0(M_1)\to E^0(F_0(M))\bigoplus E^1(M_1)\to I^0 \to 0$$ is
exact and split. So $\fd _RI^0\leq s+m$ and $\fd
_RE^i(F_0(M))\bigoplus E^{i+1}(M_1)\leq s+m+i$ for any $i\geq 1$.
Notice that $E^0(M)$ is isomorphic to a direct summand of $I^0$ and
$E^i(M)$ is isomorphic to a direct summand of $E^i(F_0(M))\bigoplus
E^{i+1}(M_1)$ for any $i\geq 1$, thus we have $\fd _RE^i(M)\leq
s+m+i$ for any $i\geq 0$. \hfill $\square$

\vspace{0.2cm}

We also need the following

\vspace{0.2cm}

{\bf Lemma 4.8.} {\it Let $M\in \Mod R$ and $n$ be a non-negative
integer.

(1) If $R$ is a right Noetherian ring and $\id _{R^{op}}F_0(M^+)\leq
\id _{R^{op}}M^++n$, then $\fd _RE^0(M)\leq \fd _RM+n$.

(2) If $R$ is a left Noetherian ring and $\id _{R^{op}}M^+\leq \id
_{R^{op}}F_0(M^+)+n$, then $\fd _RM\leq \fd _RE^0(M)+n$.}

\vspace{0.2cm}

{\it Proof.} (1) Without loss of generality, assume that $\fd
_RM=s<\infty$. Then $\id _{R^{op}}M^+=s$ by Lemma 2.4(1). So $\id
_{R^{op}}F_0(M^+)\leq \id _{R^{op}}M^+=s+n$ by assumption, and hence
we get an injective preenvelope $0\to M^{++}\to [F_0(M^+)]^+$ of
$M^{++}$ with $\fd _R[F_0(M^+)]^+=\id _{R^{op}}F_0(M^+)\leq s+n$ by
Lemma 2.4. Notice that there exists an embedding $M\hookrightarrow
M^{++}$ by [St, p.48, Exercise 41], thus $E^0(M)$ is isomorphic to a
direct summand of $[F_0(M^+)]^+$ and therefore $\fd _RE^0(M)\leq
s+n$.

(2) Without loss of generality, assume that $\fd _RE^0(M)=s<\infty$.
By Lemmas 2.3 and 2.4(1), $[E^0(M)]^+\twoheadrightarrow M^+$ is a
flat precover of $M^+$ in $\Mod R^{op}$ with $\id
_{R^{op}}[E^0(M)]^+=s$. So $F_0(M^+)$ is isomorphic to a direct
summand of $[E^0(M)]^+$ and $\id _{R^{op}}F_0(M^+)\leq s$. Then by
assumption, we have that $\id _{R^{op}}M^+\leq \id
_{R^{op}}F_0(M^+)+n\leq s+n$. It follows from Lemma 2.4(1) that $\fd
_RM\leq s+n$. \hfill$\square$

\vspace{0.2cm}

We are now in a position to state the main result in this section,
which is more general than Theorem 1.2.

\vspace{0.2cm}

{\bf Theorem 4.9.} {\it For a left Noetherian ring $R$, consider the
following conditions.

(1) $_RR$ satisfies the Auslander condition.

(2) Any flat left $R$-module satisfies the Auslander condition.

(3) $\fd _RE^i(M)\leq \fd _RM+i$ for any left $R$-module $M$ and $i
\geq 0$.

(4) $\fd _RE^0(M)\leq \fd _RM$ for any left $R$-module $M$.

(5) $\id _{R^{op}}F_i(E)\leq i$ for any injective right $R$-module
$E$ and $i \geq 0$.

(6) $\id _{R^{op}}F_i(N)\leq \id _{R^{op}}N+i$ for any right
$R$-module $N$ and $i \geq 0$.

(7) $\id _{R^{op}}F_0(N)\leq \id _{R^{op}}N$ for any right
$R$-module $N$.

We have $(1)\Leftrightarrow (2)\Leftrightarrow (3)\Leftrightarrow
(4)\Rightarrow (5) \Leftrightarrow (6) \Leftrightarrow (7)$. If $R$
is further right Noetherian, then all of the above and below
conditions are equivalent.

$(i)^{op}$ The opposite version of $(i)$ ($1\leq i \leq 7$).}

\vspace{0.2cm}

{\it Proof.} $(2)\Rightarrow (1)$ is trivial, and $(1)\Rightarrow
(2)$ follows from Corollary 3.2. $(2)\Leftrightarrow (3)
\Leftrightarrow (4)$ follow from Proposition 4.7 and Lemma 4.3(1),
and $(5)\Leftrightarrow (6) \Leftrightarrow (7)$ follow from
Proposition 4.6 and Lemma 4.3(2). By Proposition 4.5, we have
$(1)\Rightarrow (5)$.

Assume that $R$ is a left and right Noetherian ring. Then
$(1)\Leftrightarrow (1)^{op}$ follows from [FGR, Theorem 3.7], and
$(7)\Rightarrow (4)$ follows from Lemma 4.8(1). \hfill$\square$

\vspace{0.2cm}

Observe that Miyachi proved in [M, Theorem 4.1] that if $R$ is a
right coherent and left Noetherian projective $K$-algebra over a
commutative ring $K$, then $R$ satisfies the Auslander condition
(that is, $_RR$ is $G_{\infty}(0)$) if and only if $\fd _RE^0(M)\leq
\fd _RM$ for any left $R$-module $M$. Theorem 4.9 extends this
result.




By Theorems 4.9, we immediately have the following

\vspace{0.2cm}

{\bf Corollary 4.10.} {\it Let $R$ be a left Noetherian ring such
that $_RR$ satisfies the Auslander condition. If $M \in \Mod R$ with
$\fd _RM\leq s(<\infty)$, then $M$ is $G_{\infty}(s)$.}

\vspace{0.2cm}

{\it Remark 4.11.} By the dimension shifting, it is easy to verify
that the converse of Corollary 4.10 holds true when $\id _RM<\infty$
even without the assumption ``$R$ is a left and right Noetherian
ring satisfying the Auslander condition". However, this converse
does not hold true in general. For example, let $R$ be a quasi
Frobenius ring with infinite global dimension. Then $R$ is a left
and right Artinian ring satisfying the Auslander condition and every
module in $\Mod R$ is $G_{\infty}(0)$, but there exists a module in
$\Mod R$ which is not flat because $\gldim R$ is infinite.

\vspace{0.2cm}

For any $n, k\geq 0$, we use $\G_n(k)$ to denote the full
subcategory of $\Mod R$ consisting of modules being $G_n(k)$, and
denote by $\G_{\infty}(k)= \bigcap_{n\geq 0}\G_n(k)$. By [H2,
Corollary 3.9], it is easy to get the following

\vspace{0.2cm}

{\bf Proposition 4.12.} {\it Let $0\to X\to X^0\to X^1$ be an exact
sequence in $\Mod R$, and let $s\geq 0$ and $n\geq 1$. If $X^0\in
\G_n(s)$ and $X^1\in \G_{n-1}(s+1)$, then $X\in \G_n(s)$.}

\vspace{0.2cm}

For any $n\geq 0$, we use $\mathscr{F}^n(\Mod R)$ to denote the
subcategory of $\Mod R$ consisting of modules with flat dimension at most $n$.

\vspace{0.2cm}

{\bf Corollary 4.13.} {\it Let $R$ be a left Noetherian ring such
that $_RR$ satisfies the Auslander condition. Then we have

(1) $\G_{\infty}(0)=\mathscr{F}^0(\Mod R)$ if and only if $\G_{\infty}(s)=\mathscr{F}^s(\Mod R)$
for any $s\geq 0$.

(2) $\G_{\infty}(0)\bigcap\mod R=\mathscr{F}^0(\mod R)$ if and only if $\G_{\infty}(s)\bigcap\mod R
=\mathscr{F}^s(\mod R)$ for any $s\geq 0$.}

\vspace{0.2cm}

{\it Proof.} (1) The sufficiency is trivial, so it suffices to prove
the necessity. By Corollary 4.10, we have $\mathscr{F}^s(\Mod
R)\subseteq\G_{\infty}(s)$ for any $s\geq 0$. In the following we
will prove the converse inclusion by induction on $s$. The case for
$s=0$ follows from the assumption. Now suppose $s\geq 1$ and $M\in
\G_{\infty}(s)$. Let $0\to K\to F_0(M) \to M \to 0$ be an exact
sequence in $\Mod R$. By assumption $F_0(M)\in \G_{\infty}(0)$. So
$K\in \G_{\infty}(s-1)$ by Proposition 4.12, and hence $\fd_RK\leq
s-1$ by the induction hypothesis. It follows that $\fd_RM\leq s$ and
$M\in \mathscr{F}^s(\Mod R)$, which implies that
$\G_{\infty}(s)\subseteq\mathscr{F}^s(\Mod R)$.

(2) It is an immediate consequence of (1). \hfill$\square$

\vspace{0.2cm}

As applications of the results obtained above, in the rest of this
section we will study the properties of rings satisfying the
Auslander condition with finite certain homological dimension. In
particular, we will get some equivalent characterizations of
Auslander-Gorenstein rings and Auslander-regular rings.

For a module $M\in \Mod R$ and a non-negative integer $t$, we use
$\Omega ^t(M)$ to denote the $t$th syzygy of $M$ (note: $\Omega
^0(M)=M$). It is known that $\Omega ^t(M)$ is unique up to projective
equivalence for a given module $M$.

\vspace{0.2cm}

{\bf Lemma 4.14.} {\it Let $R$ be a left Noetherian ring. For a
module $M\in \Mod R$ and non-negative integers $t$ and $n$, if $\fd
_R\Omega ^t(M)\leq \fd _RE^0(\Omega ^t(M))+n$, then $\fd _RM\leq \fd
_RE^0(_RR)+n+t$.}

\vspace{0.2cm}

{\it Proof.} Let $M\in \Mod R$. Then there exist index sets $J_0,
\cdots, J_{t-1}$ such that we have the following exact sequence:
$$0 \to \Omega ^t(M) \to R^{(J_{t-1})} \to \cdots \to R^{(J_0)} \to M \to 0$$
in $\Mod R$. Because $E^0(R^{(J_{t-1})})=[E^0(_RR)]^{(J_{t-1})}$ by
[B1, Theorem 1.1] and [AF, Proposition 18.12(4)], $\fd _RE^0(R^{(J_{t-1})})
=\fd_RE^0(_RR)$. Notice that $E^0(\Omega ^t(M))$ is isomorphic to a direct
summand of $E^0(R^{(J_{t-1})})$, so $\fd _RE^0(\Omega ^t(M))\leq \fd
_RE^0(_RR)$. Thus by assumption we have that $\fd _R\Omega ^t(M)\leq
\fd _RE^0(\Omega ^t(M))+n\leq \fd _RE^0(_RR)+n$ and $\fd _RM\leq \fd
_RE^0(_RR)+n+t$. \hfill $\square$

\vspace{0.2cm}

Recall from [Bj] that a left and right Noetherian ring $R$ is called
{\it Auslander-Gorenstein} (resp. {\it Auslander-regular}) if $R$
satisfies the Auslander condition and $\id _RR=\id _{R^{op}}R$
(resp. $\gldim R)<\infty$. Also recall that $\findim R=\sup\{\pd
_RM\ |\ M\in \mod R$ with $\pd _RM<\infty\}$.

As an application of Theorem 4.9, we get some equivalent
characterizations of rings satisfying the Auslander condition with
finite left self-injective dimension as follows, which generalizes
[M, Proposition 4.4].

\vspace{0.2cm}

{\bf Theorem 4.15.} {\it For a left and right Noetherian ring $R$
and a positive integer $n$, the following statements are equivalent.

(1) $R$ satisfies the Auslander condition with $\id _RR\leq n$.

(2) $\id _{R^{op}}F_0(N)\leq \id _{R^{op}}N\leq \id
_{R^{op}}F_0(N)+n-1$ for any right $R$-module $N$ with finite
injective dimension.

(3) $\fd _RE^0(M)\leq \fd _RM\leq \fd _RE^0(M)+n-1$ for any left
$R$-module $M$ with finite flat dimension.}

\vspace{0.2cm}

{\it Proof.} $(1)\Rightarrow (2)$ Let $N\in \Mod R^{op}$ with finite
injective dimension. By Theorem 4.9, we have $\id
_{R^{op}}F_0(N)\leq \id _{R^{op}}N$. So we only need to prove the
latter inequality. Because $\id _RR\leq n$, $\id _{R^{op}}N\leq n$
by Lemma 2.5(2). So if $\id _{R^{op}}F_0(N)\geq 1$, then the
assertion holds true. Suppose $F_0(N)$ is injective. We have an
exact sequence:
$$0\to B \to F_0(N) \to N \to 0$$ in $\Mod R^{op}$ with
$\id _{R^{op}}B<\infty$. If $\id _{R^{op}}N=n$, then $\id
_{R^{op}}B=n+1$. It follows from Lemma 2.5(2) that $\id _RR\geq
n+1$, which is a contradiction. Thus we have $\id _{R^{op}}N\leq
n-1$.

$(2)\Rightarrow (3)$ Let $M\in \Mod R$ with finite flat dimension.
Then $M^+\in \Mod R^{op}$ with finite injective dimension by Lemma
2.4(1). Thus by Lemma 4.8, we get the assertion.

$(3)\Rightarrow (1)$ By (3) and Theorem 4.9, $R$ satisfies the
Auslander condition. Let $M\in \mod R$ with $\pd _RM(=\fd
_RM)<\infty$. Then $\fd _R\Omega ^1(M)<\infty$. By (3), we have $\fd
_R\Omega ^1(M)\leq \fd _RE^0(\Omega ^1(M))+n-1$. So $\pd _RM=\fd
_RM\leq \fd _RE^0(_RR)+n=n$ by Lemma 4.14. Thus we have $\findim
R\leq n$. It follows from [HI, Corollary 5.3] that $\id _RR\leq n$.
\hfill $\square$

\vspace{0.2cm}

In view of Theorem 4.15 it would be interesting to ask the following

\vspace{0.2cm}

{\bf Question 4.16.} Let $R$ be a left and right Noetherian ring
satisfying the Auslander condition with $\id _RR<\infty$. Is then
$\id _{R^{op}}R<\infty$? that is, is $R$ Auslander-Gorenstein?

\vspace{0.2cm}

By [H1, Proposition 4.6], the answer to Question 4.16 is positive if
$R$ is a left and right Artinian ring. It is a generalization of
[AR3, Corollary 5.5(b)].

Putting $n=1$ in Theorem 4.15, we have the following

\vspace{0.2cm}

{\bf Corollary 4.17.} {\it For a left and right Noetherian ring $R$,
the following statements are equivalent.

(1) $R$ satisfies the Auslander condition with $\id _RR\leq 1$.

(2) $\id _{R^{op}}F_0(N)=\id _{R^{op}}N$ for any right $R$-module
$N$ with finite injective dimension.

(3) $\fd _RE^0(M)=\fd _RM$ for any left $R$-module $M$ with finite
flat dimension.}

\vspace{0.2cm}

As another application of Theorem 4.9, we get some equivalent
characterizations of Auslander-regular rings as follows, which
generalizes [M, Corollary 4.5].

\vspace{0.2cm}

{\bf Theorem 4.18.} {\it For a left and right Noetherian ring $R$
and a positive integer $n$, the following statements are equivalent.

(1) $R$ is an Auslander-regular ring with $\gldim R\leq n$.

(2) $\id _{R^{op}}F_0(N)\leq \id _{R^{op}}N\leq \id
_{R^{op}}F_0(N)+n-1$ for any right $R$-module $N$.

(3) $\fd _RE^0(M)\leq \fd _RM\leq \fd _RE^0(M)+n-1$ for any left
$R$-module $M$.}

\vspace{0.2cm}

{\it Proof.} By Theorem 4.15 and Lemma 4.8, we have $(1)\Rightarrow
(2)\Rightarrow (3)$.

$(3)\Rightarrow (1)$ By (3) and Theorem 4.9, $R$ satisfies the
Auslander condition. Let $M\in \mod R$. By (3), we have $\fd
_R\Omega ^1(M)\leq \fd _RE^0(\Omega ^1(M))+n-1$. So $\pd _RM=\fd
_RM\leq \fd _RE^0(_RR)+n=n$ by Lemma 4.14. Thus we have $\gldim
R\leq n$. \hfill $\square$

\vspace{0.2cm}

Putting $n=1$ in Theorem 4.18, we have the following

\vspace{0.2cm}

{\bf Corollary 4.19.} {\it For a left and right Noetherian ring $R$,
the following statements are equivalent.

(1) $R$ is an Auslander-regular ring with $\gldim R\leq 1$.

(2) $\id _{R^{op}}F_0(N)=\id _{R^{op}}N$ for any right $R$-module
$N$.

(3) $\fd _RE^0(M)=\fd _RM$ for any left $R$-module $M$.}

\vspace{0.5cm}

\centerline{\large \bf 5. Approximation presentations and Gorenstein
algebras}

\vspace{0.2cm}

In this section, $R$ is an Artinian algebra. We will establish the
connection between the Auslander and Reiten's conjecture mentioned
in the introduction and the contravariant finiteness of the full
subcategory of $\mod R$ consisting of modules satisfying the
Auslander condition.

For $n\geq 0$, we use $\mathscr{I}^n(\Mod R)$ to denote the full
subcategory of $\Mod R$ consisting of modules with injective
dimension at most $n$. For a module $M\in \Mod R$, we denote by
$\Omega^{-n}(M)$ the $n$th cosyzygy of $M$. The following
approximation theorem plays a crucial role in the rest of this
section.

\vspace{0.1cm}

{\bf Theorem 5.1.} {\it Let $_RR\in\G_n(k)$ and $R_R\in\G_n(k)^{op}$
with $n,k\geq 0$. Then for any $M\in \Mod R$ and $1\leq i\leq n-1$,
there exist the following commutative diagrams with exact rows:
$$\xymatrix{0 \ar[r] & M \ar[r] \ar@{=}[d]
& I_{i+1}(M) \ar[r] \ar@{>>}[d] & G_{i+1}(M) \ar[r] \ar@{>>}[d] & 0\\
0 \ar[r] & M \ar[r] & I_{i}(M) \ar[r] & G_{i}(M) \ar[r] & 0}$$
and
$$\xymatrix{0 \ar[r] & I^{i+1}(M) \ar[r] \ar@{>>}[d]
& G^{i+1}(M) \ar[r] \ar@{>>}[d] & M \ar[r] \ar@{=}[d] & 0\\
0 \ar[r] & I^{i}(M) \ar[r] & G^{i}(M) \ar[r] & M \ar[r] & 0}$$
with $G_{j}(M), G^{j}(M)\in\G_{j}(k)$,
and $I_{j}(M), I^{j}(M)\in \mathscr{I}^{j+k}(\Mod R)$ for $j=i, i+1$.}

\vspace{0.2cm}

{\it Proof.} By [H2, Corollary 3.7 and Lemma 3.1(1)], we have the
following commutative diagrams with exact columns and rows:

{\tiny $$\xymatrix{& 0 \ar[d] & 0 \ar[d] & 0 \ar[d] & \\
0 \ar[r] & M \ar[d] \ar[r] & I_i(M) \ar[d] \ar[r] & G_i(M) \ar[d] \ar[r] & 0\\
0 \ar[r] & E^0(M) \ar[d] \ar[r] & E^0(M)\bigoplus(\bigoplus _{j=0}^{i-1}P_j(E^{j+1}(M)))
\ar[d] \ar[r] & \bigoplus _{j=0}^{i-1}P_j(E^{j+1}(M)) \ar[d] \ar[r] & 0\\
0 \ar[r] & E^1(M) \ar[d] \ar[r] & E^1(M)\bigoplus(\bigoplus _{j=0}^{i-2}P_j(E^{j+2}(M)))
\ar[d] \ar[r] & \bigoplus _{j=0}^{i-2}P_j(E^{j+2}(M)) \ar[d] \ar[r] & 0\\
& \vdots \ar[d] & \vdots \ar[d] & \vdots \ar[d] & \\
0 \ar[r] & E^{i-2}(M) \ar[d] \ar[r] & E^{i-2}(M)
\bigoplus (P_1(E^{i}(M))\bigoplus P_0(E^{i-1}(M)))
\ar[d] \ar[r] & P_1(E^{i}(M))\bigoplus P_0(E^{i-1}(M)) \ar[d] \ar[r] & 0\\
0 \ar[r] & E^{i-1}(M) \ar[d] \ar[r] & E^{i-1}(M)\bigoplus P_0(E^{i}(M))
\ar[d] \ar[r] & P_0(E^{i}(M)) \ar[d] \ar[r] & 0\\
0 \ar[r] & \Omega^{-i}(M) \ar[d] \ar[r] & E^i(M) \ar[d] \ar[r]
& \Omega^{-(i+1)}(M) \ar[d] \ar[r] & 0\\
& 0  & 0  & 0 & }$$} where
$I_i(M)=\Ker(E^0(M)\bigoplus(\bigoplus _{j=0}^{i-1}P_j(E^{j+1}(M)))
\to E^1(M)\bigoplus(\bigoplus _{j=0}^{i-2}P_j(E^{j+2}(M))))$ and $G_i(M)=\Ker
(\bigoplus _{j=0}^{i-1}P_j(E^{j+1}(M)) \to \bigoplus _{j=0}^{i-2}P_j(E^{j+2}(M)))$
for any $i\geq 1$.

Consider the following pull-back diagram:
$$\xymatrix{& & 0 \ar[d] & 0 \ar[d] & \\
0 \ar[r] & \Omega ^1(E^{i+1}(M)) \ar[r] \ar@{=}[d]& X_{i+1} \ar[r] \ar[d]
& \Omega^{-(i+1)}(M) \ar[r] \ar[d] & 0\\
0 \ar[r] & \Omega ^1(E^{i+1}(M)) \ar[r] & P_0(E^{i+1}(M)) \ar[r] \ar[d]
& E^{i+1}(M) \ar[r] \ar[d] & 0\\
& & \Omega ^{-(i+2)}(M) \ar@{=}[r] \ar[d] & \Omega ^{-(i+2)}(M) \ar[d] &\\
& & 0 & 0 &}$$ By [H2, Corollary 3.7 and Lemma 3.1(1)] again, for
any $i\geq 1$ we have the following commutative and exact columns
and rows:

{\tiny $$\xymatrix{& 0 \ar[d] & 0 \ar[d] & 0 \ar[d] & \\
0 \ar[r] & \Omega ^{i+1}(E^{i+1}(M)) \ar[d] \ar[r] & G_{i+1}(M)
\ar[d] \ar[r] & G_{i}(M) \ar[d] \ar[r] & 0\\
0 \ar[r] & P_{i}(E^{i+1}(M)) \ar[d] \ar[r] & \bigoplus _{j=0}^{i}P_j(E^{j+1}(M))
\ar[d] \ar[r] & \bigoplus _{j=0}^{i-1}P_j(E^{j+1}(M)) \ar[d] \ar[r] & 0\\
0 \ar[r] & P_{i-1}(E^{i+1}(M)) \ar[d] \ar[r] & \bigoplus _{j=0}^{i-1}P_j(E^{j+2}(M))
\ar[d] \ar[r] & \bigoplus _{j=0}^{i-2}P_j(E^{j+2}(M)) \ar[d] \ar[r] & 0\\
& \vdots \ar[d] & \vdots \ar[d] & \vdots \ar[d] & \\
0 \ar[r] & P_{2}(E^{i+1}(M)) \ar[d] \ar[r] & P_{2}(E^{i+1}(M))
\bigoplus (P_{1}(E^{i}(M))\bigoplus P_{0}(E^{i-1}(M)))
\ar[d] \ar[r] & P_{1}(E^{i}(M))\bigoplus P_{0}(E^{i-1}(M)) \ar[d] \ar[r] & 0\\
0 \ar[r] & P_{1}(E^{i+1}(M)) \ar[d] \ar[r] & P_{1}(E^{i+1}(M))\bigoplus P_0(E^{i}(M))
\ar[d] \ar[r] & P_0(E^{i}(M)) \ar[d] \ar[r] & 0\\
0 \ar[r] & \Omega^1(E^{i+1}(M)) \ar[d] \ar[r] & X_{i+1} \ar[d] \ar[r]
& \Omega^{-(i+1)}(M) \ar[d] \ar[r] & 0\\
& 0  & 0  & 0 & }$$}
Then we get the following pull-back diagram:
$$\xymatrix{& & 0 \ar[d] & 0 \ar[d]& &\\
& & \Omega ^{i+1}(E^{i+1}(M)) \ar@{=}[r] \ar[d] & \Omega ^{i+1}(E^{i+1}(M)) \ar[d]& &\\
0 \ar[r] & M \ar@{=}[d] \ar[r] & I_{i+1}(M) \ar[d] \ar[r] & G_{i+1}(M) \ar[d] \ar[r] & 0\\
0 \ar[r] & M \ar[r] & I_{i}(M) \ar[r] \ar[d] & G_{i}(M) \ar[d] \ar[r] & 0 &\\
& & 0 & 0 & & }$$ Because $R_R\in\G_n(k)^{op}$, $\id
_RP_j(E^t(M))\leq j+k$ for any $0\leq j\leq n-1$ and $t\geq 0$ by
Lemma 4.4. So from the middle column in the first diagram we get
$\id _RI_i(M)\leq i+k$ for any $1\leq i\leq n$. Because
$_RR\in\G_n(k)$, any projective module in $\mod R$ is also in
$\G_n(k)$. So by [H2, Corollary 3.9] and the exactness of the
rightmost column in the first diagram, we have $G_i(M)\in\G_i(k)$
for any $1\leq i \leq n$. Thus the above diagram is the first
desired one.

Put $I^i(M)=I_i(\Omega^1(M))$. Then we have the following push-out diagram:
$$\xymatrix{& 0 \ar[d] & 0 \ar[d]& & &\\
0 \ar[r] & \Omega^1(M) \ar[r] \ar[d] & P_0(M) \ar[d] \ar[r] & M
\ar@{=}[d] \ar[r] & 0\\
0 \ar[r] & I^i(M) \ar[r] \ar[d] & G^i(M) \ar[r] \ar[d] & M \ar[r] & 0\\
& G_i(\Omega^1(M)) \ar[d] \ar@{=}[r] & G_i(\Omega^1(M)) \ar[d] & & &\\
& 0 & 0 & & & }$$ Note that $P_0(M)\in\G_n(k)$. For any $1\leq i\leq n$,
because $G_i(\Omega^1(M))\in\G_i(k)$ by the above argument, $G^i(M)$
is also in $\G_i(k)$ by the horseshoe lemma and the exactness of the middle
column in the above diagram. By the above argument,
we have the following pull-back diagram:

$$\xymatrix{& & 0 \ar[d] & 0 \ar[d]& &\\
& & \Omega^{i+1}(E^{i+1}(\Omega^1(M))) \ar@{=}[r] \ar[d] & \Omega^{i+1}(E^{i+1}(\Omega^1(M))) \ar[d]& &\\
0 \ar[r] & P_0(M) \ar@{=}[d] \ar[r] & G^{i+1}(M) \ar[d] \ar[r] & G_{i+1}(\Omega^1(M)) \ar[d] \ar[r] & 0\\
0 \ar[r] & P_0(M) \ar[r] & G^{i}(M) \ar[d] \ar[r] & G_{i}(\Omega^1(M)) \ar[d] \ar[r] & 0\\
& & 0 & 0 & & }$$ Then the following pull-back diagram:
$$\xymatrix{& 0 \ar[d] & 0 \ar[d] & & \\
& \Omega ^{i+1}(E^{i+1}(\Omega^1(M))) \ar@{=}[r] \ar[d]
& \Omega ^{i+1}(E^{i+1}(\Omega^1(M))) \ar[d] & &\\
0 \ar[r] & I^{i+1}(M) \ar[r] \ar[d] & G^{i+1}(M) \ar[r] \ar[d] & M \ar[r] \ar@{=}[d]& 0\\
0 \ar[r] & I^{i}(M) \ar[r] \ar[d] & G^{i}(M) \ar[r] \ar[d] & M \ar[r] & 0\\
& 0 & 0 & &}$$ is the second desired one. \hfill$\square$

\vspace{0.2cm}

If $R$ satisfies the Auslander condition, then the exact sequences
$$0\to M \to I_i(M) \to G_i(M) \to 0$$ and $$0\to I^i(M)\to G^i(M)
\to M\to 0$$ in Theorem 5.1 are a left $\mathscr{I}^i(\Mod
R)$-approximation and a right $\G_i(0)$-approximation of $M$
respectively for any $1\leq i\leq n$.

\vspace{0.2cm}

{\bf Lemma 5.2.} {\it Let $X\in \mod R$ and $\{M_i\}_{i\in I}$ be a
family of left $R$-modules, where $I$ is a directed index set. Then
for any $n\geq 0$ we have
$$\Ext_R^n(\underset{i\in I}\lim M_i, X)\cong \underset{i\in I}\colim\Ext_R^n(M_i, X).$$}

\vspace{0.2cm}

{\it Proof.} Because $R$ is an Artinian algebra, any module in $\mod
R$ is pure-injective by [GT, Theorem 1.2.19]. Then the assertion
follows from [GT, Lemma 3.3.4]. \hfill$\square$

\vspace{0.2cm}

Let $M\in \Mod R$ and $n,k\geq 0$, and let
$$\cdots \to P_i(M)\to \cdots \to P_1(M) \to P_0(M)\to M \to 0$$
be a minimal projective resolution of $M$. We use $\CoG_n(k)$ to denote the full subcategory of
$\Mod R$ consisting of the modules $M$ satisfying $\id_RP_i(M)\leq i+k$ for any
$0\leq i \leq n-1$, and denote by $\CoG_{\infty}(k)=\bigcap_{n\geq 0}\CoG_n(k)$.
We use
$\mathscr{P}^n(\mod R)$ (resp. $\mathscr{I}^n(\mod R)$) to denote
the full subcategory of $\mod R$ consisting of modules with
projective (resp. injective) dimension at most $n$. As a consequence
of Theorem 5.1, we get the following

\vspace{0.2cm}

{\bf Proposition 5.3.} {\it Let $R$ satisfy the Auslander condition
and $M\in \mod R$. Then we have

(1) There exists a countably generated left $R$-module $N\in \CoG_{\infty}(0)$ and
a monomorphism $\beta: M\rightarrowtail N$ in $\Mod R$ such that $\Hom _R(\beta, T)$
is epic for any $T\in \CoG_{\infty}(0)\bigcap \mod R$.

(2) There exists a countably cogenerated left $R$-module $G\in \G_{\infty}(0)$ and
an epimorphism $\alpha: G\twoheadrightarrow M$ in $\Mod R$ such that $\Hom _R(T^{'}, \alpha)$
is epic for any $T^{'}\in \G_{\infty}(0)\bigcap \mod R$.}

\vspace{0.2cm}

{\it Proof.} (1) Let $R$ satisfy the Auslander condition. By Theorem
5.1, for any $M\in \mod R$ and $n\geq 1$, we have the following
commutative diagram with exact rows:
$$\xymatrix{0 \ar[r] & I^{n+1}(\mathbb{D}M) \ar[r] \ar@{>>}[d]
& G^{n+1}(\mathbb{D}M) \ar[r] \ar@{>>}[d] & \mathbb{D}M \ar[r] \ar@{=}[d] & 0\\
0 \ar[r] & I^{n}(\mathbb{D}M) \ar[r] & G^{n}(\mathbb{D}M) \ar[r] & \mathbb{D}M \ar[r] & 0}$$
with $G^{i}(\mathbb{D}M)\in \G_{i}(0)^{op}\bigcap \mod R^{op}$ and $I^{i}(\mathbb{D}M)\in
\mathscr{I}^{i}(\mod R^{op})$ for $i=n,n+1$, where $\mathbb{D}$ is the ordinary Matlis
duality between $\mod R$ and $\mod R^{op}$. Then we get the following commutative
diagram with exact rows:
$$\xymatrix{0 \ar[r] & M \ar[r]^{\beta_n} \ar@{=}[d] & \mathbb{D}G^{n}(\mathbb{D}M)
\ar[r] \ar@{{>}->}[d] & \mathbb{D}I^{n}(\mathbb{D}M) \ar[r] \ar@{{>}->}[d] & 0\\
0 \ar[r] & M \ar[r]^{\beta_{n+1}} & \mathbb{D}G^{n+1}(\mathbb{D}M) \ar[r]
& \mathbb{D}I^{n+1}(\mathbb{D}M) \ar[r] & 0}$$
with $\mathbb{D}G^{i}(\mathbb{D}M)\in \CoG_{i}(0)\bigcap \mod R$ and
$\mathbb{D}I^{i}(\mathbb{D}M)\in \mathscr{P}^i(\mod R)$ for $i=n,n+1$.
Put $N_n=\mathbb{D}G^{n}(\mathbb{D}M)$ and $K_n=\mathbb{D}I^{n}(\mathbb{D}M)$
for any $n\geq 1$. Then we have the following commutative diagram with exact rows:
$$\xymatrix{P_k(N_n) \ar[r] \ar[d]^{g_{n+1,n}^{k}} & P_{k-1}(N_n) \ar[r] \ar[d]^{g_{n+1,n}^{k-1}} &
\cdots \ar[r] & P_1(N_n) \ar[r] \ar[d]^{g_{n+1,n}^{1}} & P_0(N_n) \ar[r] \ar[d]^{g_{n+1,n}^{0}} &
N_n \ar[r] \ar@{{>}->}[d]^{g_{n+1,n}} & 0\\
P_k(N_{n+1}) \ar[r] & P_{k-1}(N_{n+1}) \ar[r] &
\cdots \ar[r] &  P_1(N_{n+1}) \ar[r] & P_0(N_{n+1}) \ar[r] & N_{n+1} \ar[r] & 0}$$

If $n>m$, then put
$$g_{n,m}=g_{n,n-1}g_{n-1,n-2}\cdots g_{m+1,m}$$ and
$$g_{n,m}^k=g_{n,n-1}^kg_{n-1,n-2}^k\cdots g_{m+1,m}^k.$$ In this way, for any $k\geq 0$
we get direct systems: $\{N_n, g_{n,m}\}_{n\in \mathbb{Z}^{+}}$ and
$\{P_k(N_n), g_{n,m}^k\}_{n\in \mathbb{Z}^{+}}$, where $\mathbb{Z}^{+}$ is the set
of positive integers. Because each $g_{n,m}:N_m\to N_n$ is monic, we can identity
$\underset{n\geq 1}\lim N_n$ with the direct union. It follows that $\underset{n\geq 1}\lim N_n
=\underset{n\geq t}\lim N_n$ for any $1\leq t\leq n$. Put $N=\underset{n\geq 1}\lim N_n$.
Then $N$ is countably generated.

Because $N_t\in \CoG_t(0)\bigcap \mod R$, $\id _RP_k(N_t)\leq k$ for
any $0\leq k\leq t$. So $\underset{n\geq t}\lim P_k(N_n)$ is projective and $\id _R
\underset{n\geq t}\lim P_k(N_n)\leq k$ for any $0\leq k\leq t$ by [B1, Theorem 1.1].
On the other hand, we have an exact sequence:
$$\cdots \to \underset{n\geq t}\lim P_t(N_n)\to \underset{n\geq t}\lim P_{t-1}(N_n)
\to \cdots \to \underset{n\geq t}\lim P_0(N_n)\to \underset{n\geq t}\lim N_n(=N)\to 0.$$
So $N\in \CoG_{\infty}(0)$. Put $K=\underset{n\geq t}\lim K_n$ and
$\beta=\underset{n\geq t}\lim \beta_n$. Then we get the following exact sequence:
$$0\to M \buildrel {\beta} \over \longrightarrow N \to K \to 0.$$
By Lemma 5.2, for any $T\in \CoG_{\infty}(0)\bigcap \mod R$, we have
$\Ext_R^1(K, T)\cong \Ext_R^1(\underset{n\geq t}\lim K_n, T)\cong
\underset{n\geq t}\colim\Ext_R^1(K_n, T)=0$, which implies that
$\Hom_R(\beta, T)$ is epic.

(2) Let $M\in \mod R$ and $T^{'}\in \G_{\infty}(0)\bigcap \mod R$. Then
$\mathbb{D}M\in \mod R^{op}$ and $\mathbb{D}T^{'}\in \CoG_{\infty}(0)^{op}\bigcap \mod R^{op}$.
By (1), there exists a monomorphism
$\beta: \mathbb{D}M \rightarrowtail N$ in $\Mod R^{op}$ with $N$ countably generated and
$N\in \CoG_{\infty}(0)$ such that $\Hom _{R^{op}}(\beta, \mathbb{D}T^{'})$ is epic.
Put $G=\mathbb{D}N$. Then $G$ is countably cogenerated and $\mathbb{D}\beta: G \twoheadrightarrow M
(\cong \mathbb{D}\mathbb{D}M)$ is epic in $\Mod R$ such that $\Hom_R(T^{'},\mathbb{D}\beta)
(\cong \Hom_R(\mathbb{D}\mathbb{D}T^{'},\mathbb{D}\beta))$ is also epic.
Because $N\in \CoG_{\infty}(0)^{op}$, $\id _{R^{op}}P_i(N)\leq i$ for any $i\geq 0$.
Note that $P_i(N)=\bigoplus_jP_j^i$ with all $P_j^i$ projective in $\mod R$
for any $i\geq 0$ by [Wa, Theorem 1]. So we get an exact sequence:
$$0\to G\to \prod_j\mathbb{D}P_j^0\to \prod_j\mathbb{D}P_j^1\to \cdots
\to \prod_j\mathbb{D}P_j^i \to \cdots$$ in $\Mod R$ with $\prod_j\mathbb{D}P_j^i$
injective and $\pd_R\prod_j\mathbb{D}P_j^i\leq i$ for any $i\geq 0$. It implies that
$G \in \G_{\infty}(0)$.
\hfill$\square$

\vspace{0.2cm}

Following [AR2], for a full subcategory $\mathscr{X}$ of $\mod R$ we denote by
$$\Rapp(\mathscr{X})=\{M\in \mod R\mid {\rm there\ exists\ a\ right}\
\mathscr{X}-{\rm approximation\ of}\ M\},$$
$$\Lapp(\mathscr{X})=\{M\in \mod R\mid {\rm there\ exists\ a\ left}\
\mathscr{X}-{\rm approximation\ of}\ M\}.$$
We use $\mathscr{P}^{\infty}(\mod R)$ (resp. $\mathscr{I}^{\infty}(\mod R)$) to denote
the full subcategory of $\mod R$ consisting of modules with finite projective (resp. injective)
dimension.

\vspace{0.2cm}

{\bf Proposition 5.4.} {\it Let $R$ satisfy the Auslander condition.
Then we have

(1) $\Lapp(\CoG_{\infty}(0)\bigcap\mod R)=\{M\in \mod R\mid$ there exists an exact sequence
$0\to M\to X\to Y\to 0$ with $X\in \CoG_{\infty}(0)\bigcap\mod R$ and
$Y\in \mathscr{P}^{\infty}(\mod R)\}$.

(2) $\Rapp(\G_{\infty}(0)\bigcap\mod R)=\{M\in \mod R\mid$ there exists an exact sequence
$0\to Y\to X\to M\to 0$ with $X\in \G_{\infty}(0)\bigcap\mod R$ and
$Y\in \mathscr{I}^{\infty}(\mod R)\}$.}

\vspace{0.2cm}

{\it Proof.} It is easy to see that $\Lapp(\CoG_{\infty}(0)\bigcap\mod R)
\supseteq\{M\in \mod R\mid$ there exists an exact sequence $0\to M\to X\to Y\to 0$
with $X\in \CoG_{\infty}(0)\bigcap\mod R$ and $Y\in \mathscr{P}^{\infty}(\mod R)\}$,
and $\Rapp(\G_{\infty}(0)\bigcap\mod R)\supseteq\{M\in \mod R\mid$
there exists an exact sequence $0\to Y\to X\to M\to 0$ with
$X\in \G_{\infty}(0)\bigcap\mod R$ and $Y\in \mathscr{I}^{\infty}(\mod R)\}$. So
it suffices to prove the converse inclusions.

(1) Let $M\in \Lapp(\CoG_{\infty}(0)\bigcap\mod R)$. Because $R$ satisfies the Auslander
condition, the injective cogenerator $\mathbb{D}(R_R)$ for $\Mod R$ is in
$\CoG_{\infty}(0)\bigcap\mod R$. So we may assume that $0\to M\buildrel {f} \over
\longrightarrow X^M\to Y^M\to 0$ is exact in $\mod R$ such that $f$ is a minimal
left $\CoG_{\infty}(0)\bigcap\mod R$-approximation of $M$.

Let $0\to M\buildrel {\beta} \over \to N\to K\to 0$ be an exact
sequence in $\Mod R$ as in Proposition 5.3(1) such that
$\Hom_R(\beta, T)$ is epic for any $T\in \CoG_{\infty}(0)\bigcap
\mod R$, where $N=\underset{n\geq 1}\lim N_n(=\bigcup_{n\geq 1}
N_n)$ and $K=\underset{n\geq 1}\lim K_n(=\bigcup_{n\geq 1} K_n)$.
Note that $\Hom_R(X^M,-)|_{\CoG_{\infty}(0)\bigcap\mod
R}\xrightarrow{\Hom_R(f,-)}$ \linebreak
$\Hom_R(X^M,-)|_{\CoG_{\infty}(0)\bigcap\mod R}\to 0$ is a
projective cover of $\Hom_R(X^M,-)|_{\CoG_{\infty}(0)\bigcap\mod
R}$. Because $\Hom_R(N,-)|_{\CoG_{\infty}(0)\bigcap\mod R}$ is a
projective object in the category of functors from $\Mod R$ to
Abelian groups, we have the following commutative diagram:

$$\xymatrix{\Hom_R(X^M,-)|_{\CoG_{\infty}(0)\bigcap\mod R} \ar[r]^{\Hom_R(f,-)} \ar[d]_{\Hom_R(s,-)}
& \Hom_R(M,-)|_{\CoG_{\infty}(0)\bigcap\mod R} \ar[r] \ar@{=}[d] & 0\\
\Hom_R(N,-)|_{\CoG_{\infty}(0)\bigcap\mod R} \ar[r]^{\Hom_R(\beta,-)} \ar[d]_{\Hom_R(t,-)}
& \Hom_R(M,-)|_{\CoG_{\infty}(0)\bigcap\mod R} \ar[r] \ar@{=}[d] & 0\\
\Hom_R(X^M,-)|_{\CoG_{\infty}(0)\bigcap\mod R} \ar[r]^{\Hom_R(f,-)}
& \Hom_R(M,-)|_{\CoG_{\infty}(0)\bigcap\mod R} \ar[r] & 0}$$ where
$s\in \Hom_R(N,X^M)$ and $t\in \Hom_R(X^M,N)$. Then
$\Hom_R(st,-)=\Hom_R(t,-) \Hom_R(s,-)$ is an isomorphism. So there
exist $s^{'}\in \Hom_R(K,Y^M)$ and $t^{'}\in \Hom_R(Y^M,K)$ such
that the following diagram commutes:
$$\xymatrix{0 \ar[r] & M \ar[r]^{f} \ar@{=}[d]
& X^M \ar[r] \ar[d]^{t} & Y^M \ar[r] \ar[d]^{t^{'}} & 0\\
0 \ar[r] & M \ar[r]^{\beta} \ar@{=}[d]
& N \ar[r] \ar[d]^{s} & K \ar[r] \ar[d]^{s^{'}} & 0\\
0 \ar[r] & M \ar[r]^{f} & X^M \ar[r] & Y^M \ar[r] & 0}$$
By the minimality of $f$, we have that $st$ is an isomorphism and so is $s^{'}t^{'}$.
It implies that $t^{'}:Y^M\to K(=\underset{n\geq 1}\lim K_n=\bigcup_{n\geq 1} K_n)$
is a split monomorphism. Because $Y^M$ is finitely generated, $\Im t^{'}\subseteq K_n$
for some $n$. So $Y^M$ is isomorphic to a direct summand of $K_n$ and hence
$\pd _RY^M\leq n$.

(2) Let $M\in \Rapp(\G_{\infty}(0)\bigcap\mod R)$. Then $\mathbb{D}M\in \Lapp(\CoG_{\infty}(0)^{op}
\bigcap\mod R^{op})$. By (1) there exists an exact sequence:
$$0\to \mathbb{D}M \to X\to Y\to 0$$ with $X\in \CoG_{\infty}(0)^{op}\bigcap\mod R^{op}$
and $Y\in \mathscr{P}^{\infty}(\mod R^{op})$. So we get an exact sequence:
$$0\to \mathbb{D}Y \to \mathbb{D}X \to M\to 0$$ with $\mathbb{D}X\in \G_{\infty}(0)\bigcap\mod R$
and $\mathbb{D}Y\in \mathscr{I}^{\infty}(\mod R)$.
\hfill$\square$

\vspace{0.2cm}

As a consequence of Proposition 5.4, we get the following

\vspace{0.1cm}

{\bf Proposition 5.5.} {\it Let $R$ satisfy the Auslander condition.
Then we have

(1) $\Rapp(\G_{\infty}(0)\bigcap\mod R)=\{M\in \mod R\mid$ there exists a positive integer
$n$ such that $\Omega ^{-n}(M)\in \G_{\infty}(n)\bigcap\mod R\}$.

(2) $\Lapp(\CoG_{\infty}(0)\bigcap\mod R)=\{M\in \mod R\mid$ there exists a positive
integer $n$ such that $\Omega ^{n}(M)\in \CoG_{\infty}(n)\bigcap\mod R\}$.}

\vspace{0.1cm}

{\it Proof.} (1) Let $M\in \Rapp(\G_{\infty}(0)\bigcap\mod R)$. Then
by Proposition 5.4(2), there exists an exact sequence $0\to Y\to
X\to M\to 0$ with $X\in \G_{\infty}(0)\bigcap\mod R$ and $Y\in
\mathscr{I}^{\infty}(\mod R)$. Assume that $\id _RY=k(<\infty)$.
Then for any $n>k$, $\Ext_R^1(-, \Omega^{-n+1}(X))\cong
\Ext_R^n(-,X)\cong \Ext_R^n(-, M)\cong \Ext_R^1(-,
\Omega^{-n+1}(M))$, which implies that $\Omega^{-n+1}(X)$ and
$\Omega^{-n+1}(M)$ are injectively equivalent. Because $X\in
\G_{\infty}(0)$, $\Omega^{-n+1}(X)\in \G_{\infty}(n-1)$. So
$\Omega^{-n+1}(M)\in \G_{\infty}(n-1)$ and $\Omega^{-n}(M)\in
\G_{\infty}(n)$.

Conversely, assume that $\Omega^{-n}(M)\in \G_{\infty}(n)\bigcap\mod R$.
We have the following commutative diagrams with exact columns and rows:

{\tiny $$\xymatrix{& 0 \ar[d] & 0 \ar[d] & 0 \ar[d] & \\
0 \ar[r] & I \ar[d] \ar[r] & G \ar[d] \ar[r] & M \ar[d] \ar[r] & 0\\
0 \ar[r] & K_0 \ar[d] \ar[r] & P_0(E^{0}(M))\ar[d] \ar[r] & E^{0}(M) \ar[d] \ar[r] & 0\\
0 \ar[r] & K_1 \ar[d] \ar[r] & P_0(E^{1}(M))\ar[d] \ar[r] & E^{1}(M) \ar[d] \ar[r] & 0\\
& \vdots \ar[d] & \vdots \ar[d] & \vdots \ar[d] & \\
0 \ar[r] & K_{n-2} \ar[d] \ar[r] & P_0(E^{n-2}(M))\ar[d] \ar[r] & E^{n-2}(M) \ar[d] \ar[r] & 0\\
0 \ar[r] & K_{n-1} \ar[d] \ar[r] & P_0(E^{n-1}(M))\ar[d] \ar[r] & E^{n-1}(M) \ar[d] \ar[r] & 0\\
& 0 & \Omega^{-n}(M) \ar[d] \ar@{=}[r] & \Omega^{-n}(M) \ar[d] & \\
& & 0  & 0 & }$$} where $G=\Ker(P_0(E^{0}(M))\to P_0(E^{1}(M)))$ and
$I=\Ker(K_0 \to K_1)$. Because $R$ satisfies the Auslander
condition, $P_0(E^{i}(M))$ is injective and satisfies the Auslander
condition for any $0\leq i\leq n-1$ by Theorem 4.9. So $\id_RK_i\leq
1$ for any $0\leq i\leq n-1$, and hence $\id_RI\leq n$ by the
exactness of the leftmost column in the above diagram. On the other
hand, by [H2, Corollary 3.9] and the exactness of the middle column
in the above diagram, we have that $G\in \G_{\infty}(0)\bigcap\mod
R$. Thus the exact sequence $0\to I \to G\to M\to 0$ is a right
$\G_{\infty}(0)\bigcap\mod R$-approximation of $M$ and $M\in
\Rapp(\G_{\infty}(0)\bigcap\mod R)$.




(2) It is dual to the proof of (1), so we omit it. \hfill$\square$

\vspace{0.2cm}

{\bf Corollary 5.6.} {\it Let $R$ satisfy the Auslander condition.
Then we have

(1) $\G_{\infty}(0)\bigcap\mod R$ is contravariantly finite in $\mod R$
if and only if there exists a positive integer $n$ such that
$\Omega ^{-n}(M)\in \G_{\infty}(n)\bigcap\mod R$ for any $M\in \mod R$.

(2) $\CoG_{\infty}(0)\bigcap\mod R$ is covariantly finite in $\mod R$
if and only if there exists a positive integer $n$ such that
$\Omega ^{n}(M)\in \CoG_{\infty}(n)\bigcap\mod R$ for any $M\in \mod R$.}

\vspace{0.2cm}

{\it Proof.} (1) The sufficiency follows from Proposition 5.5(1).

Conversely, let $\G_{\infty}(0)\bigcap\mod R$ be contravariantly
finite in $\mod R$ and $\{S_1, S_2, \cdots, S_t\}$ a complete set of
non-isomorphic simple $R$-modules. By Proposition 5.5(1), there
exists a positive integer $n_i$ such that $\Omega ^{-n_i}(S_i)\in
\G_{\infty}(n_i)$ for any $1\leq i\leq t$. Put $n=\max\{n_1, n_2,
\cdots, n_t\}$. Then $\Omega ^{-n}(S_i)\in \G_{\infty}(n)$ for any
$1\leq i\leq t$.

We will prove that $\Omega ^{-n}(M)\in \G_{\infty}(n)$ for any $M\in \mod R$
by induction on $\length(M)$ (the length of $M$). If $\length(M)=1$,
then $M\cong S_i$ for some $1\leq i\leq t$ and the assertion follows. Now
suppose  $\length(M)\geq 2$. Then there exists an exact sequence
$0\to S\to M \to M/S\to 0$ in $\mod R$ with $S$ simple and $\length(M/S)
<\length(M)$. By the induction hypothesis, both $S$ and $M/S$ are in
$\G_{\infty}(n)$. Then $M$ is also in $\G_{\infty}(n)$ by the horseshoe lemma.

(2) It is dual to the proof of (1), so we omit it. \hfill$\square$

\vspace{0.2cm}

Let $M\in \mod R$ and $P_1(M)\to P_0(M)\to M\to 0$ be a minimal
projective presentation of $M\in \mod R$. For any $n\geq 1$, recall
from [AB] that $M$ is called {\it $n$-torsionfree} if
$\Ext_{R^{op}}^i(\Tr M, R)=0$ for any $1\leq i\leq n$, where $\Tr
M=\Coker(P_0(M)^*\to P_1(M)^*)$ is the {\it transpose} of $M$ and
$(-)^*=\Hom_R(-,R)$. We use $\Omega^n(\mod R)$ (resp.
$\mathscr{T}_n(\mod R)$) to denote the full subcategory of $\mod R$
consisting of $n$-syzygy (resp. $n$-torsionfree) modules. Put
$\Omega^{\infty}(\mod R)=\bigcap_{n\geq 1}\Omega^n(\mod R)$ and
$\mathscr{T}_{\infty}(\mod R)=\bigcap_{n\geq 1}\mathscr{T}_n(\mod
R)$. In general, we have $\Omega^n(\mod
R)\supseteq\mathscr{T}_n(\mod R)$ for any $n\geq 1$ (cf. [AB,
Theorem 2.17]).

\vspace{0.2cm}

{\bf Lemma 5.7.} {\it If $R\in\G_{n}(0)$ with $n\geq 1$, then
$\G_{n}(0)\bigcap\mod R=\Omega^n(\mod R)=\mathscr{T}_n(\mod R)$; in
particular, if $R$ satisfies the Auslander condition, then
$\G_{\infty}(0)\bigcap\mod R=\Omega^{\infty}(\mod
R)=\mathscr{T}_{\infty}(\mod R)$.}

\vspace{0.2cm}

{\it Proof.} $\G_{n}(0)\bigcap\mod R=\Omega^n(\mod R)$ by [AR3, Proposition 5.1],
and $\Omega^n(\mod R)=\mathscr{T}_n(\mod R)$ by [AR4, Proposition 1.6 and Theorem 4.7].
\hfill$\square$

\vspace{0.2cm}


For a full subcategory $\mathscr{C}$ of $\mod R$, we denote by
$\mathscr{C}^{\bot_1}=\{M\in\mod R\mid\Ext^1_R(\mathscr{C}, M)=0\}$.

Auslander and Reiten conjectured in [AR3] that $R$ is Gorenstein
(that is, $\id_RR=\id_{R^{op}}R<\infty$) if $R$ satisfies the
Auslander condition. It remains still open. Now we are in a position
to establish the connection between this conjecture and the
contravariant finiteness of $\mathscr{G}_{\infty}(0)\bigcap\mod R$,
$\Omega^{\infty}(\mod R)$ and $\mathscr{T}_{\infty}(\mod R)$ as
follows.

\vspace{0.2cm}

{\bf Theorem 5.8.} {\it Let $R$ satisfy the Auslander condition.
Then the following statements are equivalent.

(1) $R$ is Gorenstein.

(2) $\mathscr{G}_{\infty}(0)\bigcap\mod R$ is contravariantly finite in $\mod R$.

(3) ${\rm Co}\mathscr{G}_{\infty}(0)\bigcap\mod R$ is covariantly finite in $\mod R$.

(4) $\Omega^{\infty}(\mod R)$ is contravariantly finite in $\mod R$.

(5) $\mathscr{T}_{\infty}(\mod R)$ is contravariantly finite in
$\mod R$.}


\vspace{0.2cm}

{\it Proof.} Because $R$ satisfies the Auslander condition if and
only if so does $R^{op}$, we get $(2)\Leftrightarrow (3)$. By Lemma
5.7, we have $(2)\Leftrightarrow (4)\Leftrightarrow (5)$.

$(1)\Rightarrow (2)$ Assume that $R$ is Gorenstein with
$\id_RR=\id_{R^{op}}R=n$. By [I, Proposition 1], $\pd_RE\leq n$ for
any injective left $R$-module $E$. So $\G_{\infty}(0)\bigcap\mod
R=\G_{n}(0)\bigcap\mod R$, and hence $\G_{\infty}(0)\bigcap\mod R$
is contravariantly finite in $\mod R$ by Theorem 5.1.

$(2)\Rightarrow (1)$ Assume that $\G_{\infty}(0)\bigcap\mod R$ is
contravariantly finite in $\mod R$. Then there exists a positive
integer $n$ such that $\Omega ^{-n}(M)\in \G_{\infty}(n)\bigcap\mod
R$ for any $M\in \mod R$ by Corollary 5.6, which implies that
$\G_{\infty}(0)\bigcap\mod R=\G_{n}(0)\bigcap\mod R$. Because
$\G_{n}(0)\bigcap\mod R=\mathscr{T}_n(\mod R)$ by Lemma 5.7,
$(\G_{\infty}(0)\bigcap\mod R)^{\bot_1}=(\G_{n}(0)\bigcap\mod
R)^{\bot_1} =\mathscr{T}_n(\mod R)^{\bot_1}=\mathscr{I}^n(\mod R)$
by [HI, Theorem 1.3]. On the other hand, it is easy to see that
$\mathscr{I}^{\infty}(\mod R)\subseteq (\G_{\infty}(0)\bigcap\mod
R)^{\bot_1}$. So $\mathscr{I}^{\infty}(\mod R)=\mathscr{I}^{n}(\mod
R)$ and hence $\mathscr{P}^{\infty}(\mod
R^{op})=\mathscr{P}^{n}(\mod R^{op})$. Thus $\id _{R^{op}}R\leq n$
by [HI, Corollary 5.3], which implies that $R$ is Gorenstein by
[AR3, Corollary 5.5(b)]. \hfill$\square$

\vspace{0.2cm}

As an application of Theorem 5.8, we obtain in the following result
some equivalent characterizations of Auslander-regular algebras.
Note that the converse of Corollary 4.10 does not hold true in
general by Remark 4.11. The following result also shows when this
converse holds true.

\vspace{0.2cm}

{\bf Theorem 5.9.} {\it The following statements are equivalent.

(1) $R$ is Auslander-regular.

(2) $\G_{\infty}(0)=\mathscr{P}^0(\Mod R)$.

(3) $\G_{\infty}(0)\bigcap\mod R=\mathscr{P}^0(\mod R)$.

(4) $\G_{\infty}(s)=\mathscr{P}^s(\Mod R)$ for any $s\geq 0$.

(5) $\G_{\infty}(s)\bigcap\mod R=\mathscr{P}^s(\mod R)$ for any $s\geq 0$.}

\vspace{0.2cm}

{\it Proof.} Both $(2)\Rightarrow (3)$ and $(4)\Rightarrow (5)$ are
trivial. By Corollary 4.13, we have $(2)\Leftrightarrow (4)$ and
$(3)\Leftrightarrow (5)$.

$(1)\Rightarrow (2)$ By (1) and Corollary 4.10, we have
$\mathscr{P}^0(\Mod R)\subseteq\G_{\infty}(0)$.

Assume that $\gldim R=n(<\infty)$ and $M\in \G_{\infty}(0)$. Then in a minimal injective resolution
$0\to M \to E^0(M)\to E^1(M)\to \cdots \to E^n(M)\to 0$ of $M$ in $\Mod R$, $\pd _RE^i(M)\leq i$
for any $0\leq i\leq n$. By the dimension shifting we have that $M$ is projective. It implies that
$\G_{\infty}(0)\subseteq\mathscr{P}^0(\Mod R)$.

$(5)\Rightarrow (1)$ By (5), $R$ satisfies the Auslander condition
and $\G_{\infty}(0)\bigcap\mod R=\mathscr{P}^{0}(\mod R)$ is
contravariantly finite in $\mod R$. So $R$ is Gorenstein by Theorem
5.8. Assume that $\id_{R^{op}}R=\id_RR=n(<\infty)$. Then $\pd_RE\leq
n$ for any injective left $R$-module $E$ by [I, Proposition 1]. So
for any $M\in \mod R$, $M\in \G_{\infty}(n)\bigcap\mod R$, and hence
$\pd_RM\leq n$ by (5). It follows that $\gldim R\leq n$.
\hfill$\square$

\vspace{0.5cm}

{\bf Acknowledgements.} This research was partially supported by the
Specialized Research Fund for the Doctoral Program of Higher
Education (Grant No. 20100091110034), NSFC (Grant No. 11171142) and
a Project Funded by the Priority Academic Program Development of
Jiangsu Higher Education Institutions.

\end{document}